\documentclass[12pt]{article}
\usepackage{amsfonts, amsmath, amsthm}
\usepackage[all]{xy}

\newtheorem{theorem}{Theorem}[section]
\newtheorem{lemma}[theorem]{Lemma}

\newtheorem{prop}[theorem]{Proposition}

\newcommand\AAA{\mathbb{A}}

\newcommand\QQ{\mathbb{Q}}

\newcommand\ZZ{\mathbb{Z}}

\newcommand\calO{\mathcal{O}}
\newcommand\idm{\mathfrak{m}}
\newcommand\gotho{\mathfrak{o}}

\newcommand\alg{\mathrm{alg}}
\newcommand\an{\mathrm{an}}
\newcommand\con{\mathrm{con}}

\newcommand\imm{\mathrm{imm}}
\newcommand\naive{\mathrm{naive}}

\newcommand\perf{\mathrm{perf}}

\newcommand\beq{\begin{equation}}
\newcommand\eeq{\end{equation}}

\newcommand\fp{\frac{1}{p}}

\newcommand\GK{\Gamma^K}
\newcommand\GL{\Gamma^L}
\newcommand\Gperf{\Gamma^{\perf}}

\newcommand\Gimm{\Gamma^{\imm}}
\newcommand\Galg{\Gamma^{\alg}}

\newcommand\Galgancon{\Galg_{\an, \con}}
\newcommand\Galgcon{\Galg_{\con}}

\newcommand\Gancon{\Gamma_{\an,\con}}
\newcommand\Gcon{\Gamma_{\con}}

\newcommand\Gimmcon{\Gimm_{\con}}
\newcommand\GKancon{\GK_{\an,\con}}
\newcommand\GKcon{\GK_{\con}}
\newcommand\GLancon{\GL_{\an,\con}}
\newcommand\GLcon{\GL_{\con}}
\newcommand\Gperfcon{\Gperf_{\con}}

\newcommand\be{\mathbf{e}}

\newcommand\bv{\mathbf{v}}
\newcommand\bw{\mathbf{w}}

\newcommand\by{\mathbf{y}}
\newcommand\bz{\mathbf{z}}
\newcommand\calE{\mathcal{E}}

\DeclareMathOperator{\Gal}{Gal}
\DeclareMathOperator{\Hom}{Hom}

\DeclareMathOperator{\FaIsoc}{\mathit{F^a}-Isoc}
\DeclareMathOperator{\FaIsocd}{\mathit{F^a}-Isoc^{\dagger}}
\DeclareMathOperator{\rank}{rank}

\begin{document}

\title{Full faithfulness for overconvergent $F$-isocrystals}
\author{Kiran S. Kedlaya \\ University of California, Berkeley}
\date{June 1, 2003}

\maketitle

\begin{abstract}
Let $X$ be a smooth variety over a field of characteristic $p>0$.
We prove that the forgetful functor from the category of overconvergent
$F$-isocrystals on $X$ to that of convergent $F$-isocrystals is
fully faithful.
The argument uses the quasi-unipotence theorem
for overconvergent $F$-isocrystals (recently proved independently by
Andr\'e, Mebkhout, and the author),
plus arguments of de~Jong. In the process, we establish a theorem
of Quillen-Suslin type (i.e., every finite projective module is free)
over rings of overconvergent power series.
\end{abstract}

\section{Introduction}

Isocrystals, and more specifically $F$-isocrystals, were constructed to
serve as $p$-adic analogues of local systems in the complex topology and
lisse $l$-adic sheaves in the \'etale topology. Unfortunately, the category
of convergent isocrystals suffers from certain pathologies, arising from
the fact that the integral of a $p$-adic analytic function on a closed
disc is no longer bounded near the boundary of the disc. For example, computed
using convergent isocrystals, the cohomology of the affine line is infinite
dimensional. 

As first noted in the work of Dwork, pointed out more
explicitly by Monsky and Washnitzer in the study of dagger cohomology,
and systematized by Berthelot in his theory of rigid cohomology, one needs
to work instead with objects satisfying an ``overconvergence'' condition.
Adding this condition eliminates the pathologies mentioned above; in fact,
Berthelot \cite{bib:ber1} showed that the rigid cohomology with constant
coefficients of an arbitrary variety over a field of characteristic $p$
is finite dimensional. (Finite dimensionality
has now also been shown for rigid cohomology with coefficients
in an $F$-isocrystal, which itself must be overconvergent \cite{bib:mefinite}.)
However, additional complications arise from the fact
that it is not always clear how to ``descend'' certain constructions from
the convergent category to the smaller overconvergent category.

The main result of this paper is a descent theorem in this vein,
which resolves a conjecture formulated by Tsuzuki
\cite[Conjecture~1.2.1(TF)]{bib:tsu4}.
Let $k$ be a field of characteristic $p>0$, $K$ the fraction field of
a complete mixed characteristic
discrete valuation ring with residue field $k$, and $X$ a smooth separated
$k$-scheme of finite type over $K$. Then for each integer $a>0$,
we can construct the
categories $\FaIsoc(X/K)$ and $\FaIsocd(X/K)$ of convergent and overconvergent
$F^a$-isocrystals, respectively, on $X$; see for example Berthelot 
\cite{bib:ber2}. One consequence of the construction is that
there is a natural forgetful functor $j^*:
\FaIsocd(X/K) \to \FaIsoc(X/K)$.
\begin{theorem} \label{thm:main}
The forgetful functor $j^*: \FaIsocd(X/K) \to \FaIsoc(X/K)$ is fully faithful.
\end{theorem}

In the case of unit-root $F$-isocrystals, this was proved for $X$ admitting
a smooth compactification by Tsuzuki \cite{bib:tsu4} and for general $X$
by Etesse \cite{bib:et}.
Tsuzuki further conjectures that the same statement
holds for isocrystals without Frobenius structure, but the methods of this
paper do not extend to that case.

The aforementioned results of Tsuzuki and Etesse have the effect of reducing
the unit-root case of Theorem~\ref{thm:main} to a local assertion about
unit-root $F$-isocrystals, which had earlier been proved by Tsuzuki
\cite{bib:tsu2}. We prove Theorem~\ref{thm:main} by using a reduction of a
similar spirit to bring the global problem down to
an analogous local assertion about $F$-isocrystals. Our proof of said
local assertion closely follows
de~Jong's proof of the equal characteristic analogue
of Tate's extension theorem for $p$-divisible groups. The key new ingredient
is the proof of Crew's conjecture \cite{bib:crew2} that overconvergent
$F$-isocrystals are quasi-unipotent. This conjecture has recently
been proved independently
by Andr\'e \cite{bib:andre}, Mebkhout \cite{bib:mebkhout}
and the author \cite{bib:me7}.

One other result of this paper may be of independent interest: we
prove (Theorem~\ref{thm:qs})
an analogue of the Quillen-Suslin theorem over any ring
of overconvergent power series over a complete discrete valuation
ring over field. That is, over such a ring, every finite
projective module is free.

\subsection*{Acknowledgments}
This work is partially based on the author's doctoral dissertation 
\cite{bib:methesis}, written under Johan de~Jong. The author was
supported by a National Science Foundation Postdoctoral Fellowship.
He also thanks the organizers of the Dwork Trimester for
their hospitality, Jean-Yves Etesse for providing a copy of his
preprint \cite{bib:et}, and Laurent Berger and Bernard le Stum
for helpful discussions.

\section{Definitions and notations}

This paper is a companion paper to \cite{bib:me7}, and so we
will adopt its notation and terminology.
For the convenience of the reader, we recall the particular definitions
and notations we will need; however, we will defer to
\cite[Sections~2--3]{bib:me7} for the verification of various compatibilities.

For $k$ a field of characteristic $p>0$, let $C(k)$ be a Cohen ring for
$k$, that is, a complete discrete valuation ring with
fraction field of characteristic 0, maximal ideal generated by $p$, and
residue field isomorphic to $k$. (The ring $C(k)$ is unique up to
noncanonical isomorphism if $k$ is not perfect; if $k$ is perfect,
$C(k)$ is canonically isomorphic to the ring $W(k)$ of Witt vectors
over $k$.)
Let $\calO$ be a finite totally ramified extension of $C(k)$,
let $\pi$ be a uniformizer of $\calO$, and
fix a ring endomorphism $\sigma_0$
on $\calO$ lifting the absolute Frobenius $x \mapsto x^p$ on $k$.
Let $q = p^f$ be a power of $p$ and put $\sigma = \sigma_0^f$.
Let $v_p$ denote the valuation on $\calO$ normalized so that
$v_p(p) = 1$, and let $|\cdot|$ denote the norm on $\calO$
given by $|x|= p^{-v_p(x)}$. For $x \in \calO$, let
$\overline{x}$ denote its reduction in $k$.

Let $\Gamma$ denote the ring of bidirectional power series
\[
\left\{\sum_{i \in \ZZ} x_i u^i: x_i \in \calO, \qquad
\lim_{i \to -\infty} v_p(x_i)
= \infty\right\}.
\]
Then $\Gamma$ is a complete discrete valuation ring, whose residue field
we identify with $k((t))$ by identifying the reduction of $\sum x_i u^i$
with $\sum \overline{x_i} t^i$.
Choose an extension of $\sigma_0$, as defined
on $\calO$, to a ring endomorphism of $\Gamma^K$ lifting the $p$-th
power $x \to x^p$ on $K$, and mapping the subring
\[
\Gcon = \left\{\sum_{i \in \ZZ} x_i u^i: x_i \in \calO, \qquad
\liminf_{i \to -\infty}
\frac{v_p(x_i)}{-i} > 0 \right\}
\]
into itself; it suffices to check that $u^{\sigma_0} \in \Gcon$.
Again put $\sigma = \sigma_0^f$.

If $k$ is perfect, we can define a functor from the category of
algebraic (finite or infinite) extensions of $k((t))$ to
the category of complete discrete valuation rings which are unramified
over $\calO$. Let $\Gamma^K$ denote the image of an extension $K$;
if $K = k((t))$ we will omit the superscript entirely, while
if $K = k((t))^{\perf}$ or $K = k((t))^{\alg}$,
we will call the image $\Gamma^{\perf}$ or $\Galg$, respectively.
If $K$ is finite over $k((t))$, then $K$ is itself isomorphic to
$k((t'))$ for some $t'$, and analogously the ring $\Gamma^K$ is 
abstractly isomorphic to $\Gamma$.

Each ring $\Gamma^K$ comes equipped with a canonical extension of
$\sigma$.
In fact, the choice of the functor depends on $\sigma_0$, at least
for $K$ which are not separable over $k((t))$: for example,
$\Gamma^{\perf}$ is the direct limit of $\Gamma \to \Gamma \to \cdots$,
where the transition maps are equal to $\sigma$.
The functoriality is more limited if $k$ is not perfect, but 
for our purposes, it will suffice to
embed $\Gamma$ into $\Gamma^K$ for $K = k^{\alg}((t))$, and then
construct $\Gamma^L$ functorially for extensions $L$ of $K$.

We define the partial valuations on
$\Galg[\fp]$ as follows: for
$x \in \Galg[\fp]$, write $x = \sum_{i=0}^\infty \pi^i [u_i]$, where the
brackets denote Teichm\"uller lifts, and put
\[
v_n(x) = \min_{v_p(\pi^i) \leq n} \{ v_t(u_i) \},
\]
where the bar denotes reduction to $k((t))^{\alg}$
and $v_t$ denotes the canonical extension of the valuation on
$k((t))$, normalized so that $v_t(t) = 1$. 
For each $r>0$, let $\Gamma^K_r$ denote the subring
of $x \in \Gamma^K$ such that
$\lim_{n \to \infty} (r v_n(x) + n) = \infty$,
and define the valuation
\[
w_r(x) = \min_n \{r v_n(x) + n\}
\]
on $\Gamma^K_r$. Define $\Gamma^K_{\con} = \cup_{r>0} \Gamma^K_r$.
It can be shown that for $K=k((t))$, the ring $\Gamma^K_{\con}$
coincides with the ring $\Gcon$ defined above. More precisely,
for $r$ sufficiently small (depending on the choice of $\sigma_0$),
the $w_r$ coincide with their ``na\"\i ve'' counterparts, in whose definition
$v_n$ is replaced by
\[
v_n^{\naive}(x) = \min_{j: v_p(x_j) \leq n} \{j\},
\]
where $x = \sum_j x_j u^j$.

We can define a Fr\'echet topology on $\Gamma^K_r[\fp]$ using
the norms $w_s$ for $0 < s \leq r$, and take $\Gamma^K_{\an, r}$
as the Fr\'echet completion of $\Gamma^K_r[\fp]$. (However, this is only
known to behave well when $K$ is either finite over $k((t))$ or perfect.)
Put $\Gamma^K_{\an,\con} = \cup_{r>0} \Gamma^K_{\an,r}$.
For $k=k((t))$, the ring $\Gamma^K_{\an,\con} = \Gancon$
coincide with the Robba ring,
the ring of germs of functions analytic on some open $p$-adic annulus
with outer radius 1. Concretely, this ring consists of series
\[
\left\{ \sum_j x_j u^j: x_j \in \calO[1/p], 
\qquad \liminf_{j \to -\infty}
\frac{v_p(x_j)}{-j} > 0, \qquad \liminf_{j \to \infty} 
\frac{v_p(x_j)}{j} \geq 0\right\}.
\]
It turns out that $\GKancon$ is a B\'ezout ring (every finitely generated
ideal is principal) if $K$ is finite over $k((t))$ or $K$ is perfect,
by \cite[Theorem~3.12]{bib:me7}, but we will not explicitly need
this fact.

We define a \emph{$\sigma$-module} over an integral domain $R$
to be a finite free
$R$-module $M$ equipped with an $R$-linear map $F: M \otimes_{R, \sigma} R
\to M$ that becomes an isomorphism over $R[\fp]$; the tensor product
notation indicates that $R$ is viewed as an $R$-module via $\sigma$.
To specify $F$, it is equivalent to specify an additive,
$\sigma$-linear map from $M$ to $M$ that acts on any basis of $M$ by
a matrix invertible over $R[\fp]$. We abuse notation and refer to this map
as $F$ as well; since we will only use the $\sigma$-linear map in what
follows, there should not be any confusion induced by this.

For $K = k((t))$ and $R$ one of $\GK$, $\GK[\fp]$,
$\GKcon$, $\GKcon[\fp]$, $\GKancon$, let $\Omega^1_R$ be the free module
over $R$ generated by a symbol $du$, and define the derivation
$d: R \to \Omega^1_R$ by the formula
\[
d\left( \sum_j x_j u^j \right) = \left( \sum_j j x_j u^{j-1} \right)\,du.
\]
Likewise, if $K$ is a finite extension of $k((t))$, we can make the
same definition by writing elements as power series in terms of some
$u' \in \GKcon$ which lifts a uniformizer of $K$.

We define a $(\sigma, \nabla)$-module over $R$
to be a $\sigma$-module
$M$ plus a connection $\nabla: M \to M \otimes_R \Omega^1_{R}$ 
(i.e., an additive map satisfying the Leibniz rule
$\nabla(c\bv) = c\nabla(\bv) + \bv \otimes dc$) that makes the following 
diagram commute:
\[
\xymatrix{
M \ar^-{\nabla}[r] \ar^{F}[d] & M \otimes \Omega^1_{R} \ar^{F \otimes d\sigma}[d] \\
M  \ar^-{\nabla}[r] & M \otimes \Omega^1_{R}
}
\]
Given a $\sigma$-module or $(\sigma, \nabla)$-module $M$ and an
integer $\ell$, we define the \emph{Tate twist} $M(\ell)$ as the module $M$
with the action of $F$ multiplied by $q^{\ell}$. Note that the dual
$M^* = \Hom(M, R)$ does not generally have the structure of
a $(\sigma, \nabla)$-module over $R$ (only over $R[\fp]$),
but its Tate twist $M^*(\ell) = \Hom(M, R(\ell))$ does for some $\ell$.
If $\bv$ is an element of $M$ such that 
$F\bv = \lambda \bv$ for some $\lambda \in \calO$, we say
$\bv$ is an \emph{eigenvector} of $M$ of eigenvalue $\lambda$
and slope $v_p(\lambda)$. 

There are two ways to associate a Newton polygon to a $\sigma$-module.
For $R = \Gamma$ or $R = \Gcon$, the Dieudonn\'e-Manin classification
states that
a $\sigma$-module over $R$ acquires a basis of eigenvectors over
$\Gamma^{\alg}$; the slopes of these eigenvectors are called the
\emph{generic slopes} of $M$. For $R = \Gcon$ or $R = \Gancon$, 
\cite[Theorem~3.12]{bib:me7}
states that a $\sigma$-module over $R$ acquires a basis of eigenvectors
over $\Galgancon$; the slopes of these eigenvectors are called the
\emph{special slopes} of $M$.

We use the following refinement of the Dieudonn\'e-Manin
classification, due originally to de~Jong \cite[Proposition~5.5]{bib:dej1} and
appearing also as \cite[Proposition~5.7]{bib:me7}.
\begin{prop}[Descending slope filtration] \label{prop:descfilt}
Let $M$ be a $\sigma$-module over $\Gcon$, for $k$ algebraically closed.
Then there exists a finite extension $\calO'$ of $\calO$ such that
over $\Galgcon \otimes_{\calO} \calO'$, $M$ admits a basis $\bv_1,
\dots, \bv_n$ such that $F\bv_i = \lambda_i \bv_i + \sum_{j<i} c_{ij} \bv_j$
for some $\lambda_i \in \calO$ and $c_{ij} \in \Galgcon$, with
$v_p(\lambda_1), \dots, v_p(\lambda_n)$ equal to the sequence of
generic slopes of $M$ in descending order. Moreover, if
$v_p(\lambda_1) = \cdots =v_p(\lambda_m)$, we can take $c_{ij} = 0$
for $i \leq m$.
\end{prop}

We can also define $\sigma$-modules and $(\sigma, \nabla)$-modules over
rings $R$ which are not discrete valuation rings. In particular, one typically
makes these definitions for $R$ a (smooth) \emph{dagger algebra}, a quotient of
a ring $\calO\langle x_1, \dots, x_n \rangle^\dagger$ of overconvergent
power series which is smooth over $\calO$,
in which each series $\sum_I c_I x^I$ satisfies
$\liminf v_p(c_I)/|I| > 0$. (Here $|I|$ denotes the sum of the indices
in the index set.) One defines a $(\sigma, \nabla)$-module over $R$
as above, modulo the following changes:
\begin{enumerate}
\item The underlying module $M$ is allowed to be locally free, but
not necessarily free.
\item
The module of differentials $\Omega^1_{R}$ is now the $R$-algebra obtained
from the free module over $\calO\langle x_1, \dots, x_n \rangle^\dagger$
generated by $dx_1, \dots, dx_n$ by quotienting by the necessary relations.
(The result does not depend on the presentation of $R$.) 
\item
The connection must now satisfy an integrability condition: if
$\nabla_1: M \otimes \Omega^1_R \to M \otimes \wedge^2 \Omega^1_R$ is the natural map
induced by $\nabla$, then $\nabla_1 \circ \nabla = 0$.
\end{enumerate}

\section{Splitting an exact sequence}

We first recall \cite[Proposition~3.19(c)]{bib:me7}.
\begin{prop} \label{prop:eq}
For $\lambda \in \calO$ not a unit and $x \in \GKancon$, there
is at most one $y \in \GKancon$ such that $\lambda y^\sigma - y = x$.
Moreover, if $x \in \GKcon$, then so is $y$.
\end{prop}
Our next lemma generalizes this lemma, using the descending slope filtration.
One can also give a direct proof within $\Gancon$;
we leave this as an exercise to the interested reader.
\begin{lemma} \label{lem:inj}
Let $M$ be a $\sigma$-module over $\GKcon[\fp]$, for $K$ a finite
extension of $k((t))$, whose generic
slopes are all positive.
Then the map $\bv \mapsto F\bv - \bv$ on $M \otimes_{\GKcon[1/p]} \GKancon$ is
injective. Moreover, if $F \bv - \bv \in M$, then $\bv \in M$.
\end{lemma}
\begin{proof}
For both assertions, it suffices to assume $k$ is algebraically closed.
We prove the second assertion first.
Let $\bv_1, \dots, \bv_n$ be the basis of $M$ over
$\Galgcon \otimes_{\calO} \calO$'
given by Proposition~\ref{prop:descfilt}, with $F\bv_i = \lambda_i \bv_i + \sum_{j<i} 
c_{ij} \bv_j$. Since $v_p(\lambda_i)$ are the generic slopes of $M$,
they are all positive by hypothesis.
Write $\bv = \sum_i e_i \bv_i$ with $e_i \in \Galgancon$
and put $F\bv -\bv = \sum_i
f_i \bv_i$, with $f_i \in \Galgcon$. Then
\begin{align*}
f_n &= \lambda_n e_n^\sigma - e_n \\
f_{n-1} &= \lambda_{n-1} e_{n-1}^\sigma - e_{n-1} + c_{n(n-1)} e_n^\sigma \\
&\vdots \\
f_1 &= \lambda_1 e_1^\sigma - e_1 + c_{21} e_2^\sigma + \cdots +
c_{n1} e_n^\sigma.
\end{align*}
By Proposition~\ref{prop:eq}, the map 
$x \mapsto \lambda_i x^\sigma - x$
on $\Galgancon \otimes_{\calO} \calO'$
is injective, and if $\lambda_i x^\sigma - x \in \Galgcon[\fp]$, 
then $x \in \Galgcon[\fp]$.
Applying this fact repeatedly to the above equations, we deduce successively
that $e_n, e_{n-1}, \dots, e_1$ all lie in $\Galgcon[\fp]$. 
Therefore $\bv$ is defined
over $\Galgcon[\fp] \cap \Gancon = \Gcon[\fp]$.

To establish the first assertion, suppose $F\bv = \bv$;
we then have the same equations as above, but with $f_i = 0$ for all $i$.
By Proposition~\ref{prop:eq} again, we have $e_n = 0$, then $e_{n-1} =0$,
and so on. Thus $\bv = 0$, as desired.
\end{proof}

We now use the previous lemma, plus the quasi-unipotence theorem,
to prove the following proposition, which is analogous to
\cite[Proposition~7.1]{bib:dej1}. (In de~Jong's argument, the role of
the quasi-unipotence theorem is played by Dwork's trick.)
\begin{prop} \label{prop:split}
  Let $M$ be a $(\sigma, \nabla)$-module over $\Gcon$.
Suppose there exists an exact sequence $0 \to M_1 \to M \to M_2 \to 0$
of $(\sigma, \nabla)$-modules over $\Gcon$ such that the generic slopes
of $M_1$ are all greater than the generic slopes of $M_2$. Then
the exact sequence splits over $\Gcon[\fp]$.
\end{prop}
\begin{proof}
Let $n_1$ and $n_2$ be the ranks of $M_1$ and $M_2$, respectively.
Choose a basis of $M_1$ over $\Gcon[\fp]$
and extend it to a basis of $M$. Then $\sigma$
acts on this basis via a block matrix of the form
$\begin{pmatrix} A & B \\ 0 & D \end{pmatrix}$.
To show that the exact sequence splits over a ring $R$ containing
$\Gcon[\fp]$, it is necessary and sufficient
to find a matrix $X$ over $R$ such that
\begin{equation}
  \label{eq:sigma}
-X + A X^\sigma D^{-1} = B,
\end{equation}
as then we have
\[
\begin{pmatrix}
  I_{n_1} & -X \\ 0 & I_{n_2}
\end{pmatrix} 
\begin{pmatrix}
  A & B \\ 0 & D
\end{pmatrix}
\begin{pmatrix}
  I_{n_1} & X^\sigma \\ 0 & I_{n_2}
\end{pmatrix}
= 
\begin{pmatrix}
  A & 0 \\ 0 & D
\end{pmatrix}.
\]

We first solve (\ref{eq:sigma}) over a somewhat larger ring than
$\Gcon[\fp]$,
using the quasi-unipotence theorem. For any finite extension $L$ of $K
= k^{\alg}((t))$
and any lift $u \in \GLcon$ of a uniformizer of $L$,
define $\GLancon[\log u]$ as the polynomial ring over $\GLancon$ in
an indeterminate called ``$\log u$''. Extend $\sigma$ from $\GLancon$ to
$\GLancon[\log u]$ by setting
\[
(\log u)^\sigma = p \log u + \log (u^\sigma/u^p),
\]
using the power series expansion of $\log(1+x)$ to define the second factor,
and extend the derivation $d$ to $\log u$ by setting
$d(\log u) = \frac{1}{u} \otimes du$.
By the quasi-unipotence theorem (as formulated in 
\cite[Theorem~6.13]{bib:me7}), there exists a finite separable extension
$L$ of $K$, a finite extension $\calO'$ of $\calO$, and a
basis $\bv_1, \dots, \bv_{n_1}, \bw_1, \dots, \bw_{n_2}$
of $M$ over $R_1 = \GLancon[\log u] \otimes_{\calO} \calO'$ such that
$\bv_1, \dots, \bv_{n_1}$ form a basis of $M_1$, and
\begin{align*}
  F\bv_i &= \lambda_i \bv_i \\
\nabla\bv_i &= 0 \\
F\bw_i &= \mu_i \bw_i + \sum_{j=1}^m W_{ij} \bv_j  \\
\nabla \bw_i &= \sum_{j=1}^{n_1} Y_{ij} \bv_j \otimes du
\end{align*}
for some $\lambda_1, \dots, \lambda_{n_1}, \mu_1, \dots, \mu_{n_2}$ in
$\calO$ and $W_{ij}, Y_{ij} \in R_1$.

The ring $R_1$ has the property that the map $\frac{d}{du}: R_1 \to R_1$
is surjective. (The point is that every expression $u^i (\log u)^j$
has an antiderivative, by integration by parts.) Thus we can choose $Z_{ij}
\in R_1$ such that $\frac{d}{du} Z_{ij} = Y_{ij}$. Put
$\by_i = \bw_i - \sum_{j=1}^{n_1} Z_{ij} \bv_j$.
Now $\nabla \by_i = 0$ for $i=1, \dots, n$, and we may write
$F\by_i = \mu_i \by_i + \sum_{j=1}^{n_1} V_{ij} \bv_j$ for some
$V_{ij} \in R_1$. By the compatibility relation between
$F$ and $\nabla$, we have $0 =
\nabla F\by_i = \sum_{j=1}^{n_1} \bv_j \otimes dV_{ij}$, so $V_{ij} \in \calO$
for all $i,j$. Choose $U_{ij} \in \calO$ such that $\mu_i U_{ij}^\sigma - 
\lambda_j U_{ij}
= V_{ij}$, and put $\bz_i = \by_i - \sum_{j=1}^{n_1} U_{ij} \bv_j$.
Then $F\bz_i = \mu_i \bz_i + \sum_{j=1}^{n_1} (V_{ij} + \mu_i U_{ij} -
\lambda_j U_{ij}^\sigma) = \mu_i \bz_i$, so $M$ splits over
$R$ as $M_1$
plus the span of $\bz_1, \dots, \bz_{n_2}$. 
Thus (\ref{eq:sigma}) admits a solution $X_1$ over $R_1$.

We now descend the solution $X_1$ we just found to smaller and smaller
rings. Before doing so, note that we can give the set of $n_1 \times n_2$
matrices $X$ the structure of a $\sigma$-module with Frobenius
given by $X \mapsto AX^\sigma D^{-1}$. Each generic slope of this
$\sigma$-module is a generic slope of $A$ minus a generic slope of $D$,
and thus is positive. We thus can apply Lemma~\ref{lem:inj} to this
$\sigma$-module.

Now for the descent.
First put $R_2 = \GLancon \otimes_{\calO} \calO'$,
so that $R_1 = R_2[\log u]$.
If $X \neq 0$, we can write $X = \sum_{i=0}^M Y_i (\log u)^i$ with 
$Y_M \neq 0$. If $M > 0$, we have $-Y_M + p^m A Y_M^\sigma D^{-1} = 0$;
by Lemma~\ref{lem:inj}, this forces $Y_M = 0$, contradiction.
Thus we must have $M = 0$, so $X$ is defined over $R_2$.

Next, put $R_3 = \GLancon$, so that $R_2 = R_3 \otimes_{\calO} \calO'$.
Choose a basis $c_1, \dots, c_r$ of $\calO'$ over $\calO$
with $c_1 = 1$. Then we can write $X_2$ as
a linear combination of $c_1, \dots, c_r$ whose coefficients are
matrices over $R_1$. The coefficient of $c_1$ must then also 
be a solution of \eqref{eq:sigma}; since $X$ is unique by
Lemma~\ref{lem:inj}, it must coincide with its coefficient of $c_1$.
That is, $X$ is defined over $R_3$.

Next, put $R_4 = \GLcon[\fp]$. We have $- X + A X^\sigma D^{-1} = B$;
by Lemma~\ref{lem:inj}, since $B \in \GLcon[\fp] = R_4$, we must have $X \in
R_4$.

Next, put $R_5 = \GKcon[\fp]$. Now $\Gal(L/K)$ acts on $\GLcon$ with 
fixed ring $\GKcon$. In particular, $\Gal(L/K)$ acts on the set of
solutions of \eqref{eq:sigma}. But by Lemma~\ref{lem:inj}, there is
only one solution $X$. Thus $X$ is fixed by $\Gal(L/K)$ and so is
defined over $R_5$.

Finally, put $R_6 = \Gcon[\fp]$. We now have a $\sigma$-module $Y$ over
$\Gcon[\fp]$ (the space of $n_1 \times n_2$ matrices) with all
generic slopes positive, and an eigenvector $X$ over $\GKcon[\fp]$.
By \cite[Proposition~5.3]{bib:me7}, $Y$ admits a basis $Z_1,
\dots, Z_{n_1n_2}$ on which Frobenius
acts by a matrix over $\Gcon$ with positive valuation. Write
$X = \sum c_i Z_i$, suppose there exists an integer $m$ such that
$c_i$ is not congruent to an element of $\Gcon[\fp]$ modulo $\pi^m$
for some $i$, and choose the smallest
such $m$. Write each $c_i = e_i + f_i$ with $e_i \in \Gcon[\fp]$
and $f_i \not \equiv 0 \pmod{\pi^m}$ for some $i$, and put
$X_1 = \sum f_i Z_i$.
Then $X_1 - AX_1^{\sigma}D^{-1} = \sum d_i Z_i$ with
$d_i \in \Gcon[\fp]$, but $d_i - f_i \equiv 0 \pmod{\pi^m}$
because the matrix on which $F$ acts on the $Z_i$ has positive
valuation. Then $c_i \equiv e_i + d_i \pmod{\pi^m}$, contrary to the
choice of $m$. We conclude that the $c_i$ are congruent to elements
of $\Gcon[\fp]$ modulo every power of $\pi$, so belong to $\Gcon[\fp]$.
In other words, $X$ is defined over $R_6$.

To conclude, we have shown that \eqref{eq:sigma} admits a solution
over $\Gcon[\fp]$. Thus the exact sequence splits, as desired.
\end{proof}

\section{Equality of kernels}

This section is nearly a carbon copy of \cite[Section~8]{bib:dej1},
to the extent that we have reproduced its title. The main changes are that
we work with an arbitrary $\sigma$ rather than a ``standard'' $\sigma$,
which sends some $u \in \GKcon$ to $u^p$, and that we expose
\cite[Proposition~8.1]{bib:dej1} using the technical device of
``generalized power series'', which we hope provides a small clarification
of the argument.

For $K$ a valued field
over an algebraically closed field of characteristic $p>0$, the maximal
immediate extension $K^{\imm}$, in the sense of Kaplansky \cite{bib:kap},
is the maximal extension of $K$ with value group $\QQ$.
Kaplansky shows this field can
be identified with Hahn's field of generalized power series
$x = \sum_i x_i t^i$, where $x_i \in k$ for each $i \in \QQ$, 
and the set $I = I(x)$
of indices $i$ such that $x_i \neq 0$ is well-ordered for each $x$.
The corresponding ring $\Gimm$ can be described as the ring of generalized
power series $y = \sum_i y_i w^i$, where $y_i \in \calO$ for each
$i \in \QQ$ and for each $n$, the set $I = I(n,x)$ such that
$v_p(y_i) \leq n$ is well-ordered; the partial valuation $v_n(y)$ is equal
to the smallest $i$ such that $v_p(y_i) \leq n$. Note that
$w$ has $p^n$-th roots for all $n$ and so must be the Teichm\"uller lift
of its reduction.

The following proposition corresponds to \cite[Proposition~8.1]{bib:dej1}.
\begin{prop}[after de~Jong]
For $k$ algebraically closed,
the multiplication map $\Galgcon \otimes_{\Gcon} \Gamma
\to \Galg$ is injective.  
\end{prop}
\begin{proof}
Suppose, by way of contradiction, that
$\sum_{i=1}^n f_i \otimes g_i$
is a nonzero element of $\Galgcon \otimes_{\Gcon}
\Gamma$ such that $\sum f_ig_i$ in $\Galg$, and assume $n$ is minimal for 
the existence of such an element. Then the $g_i$ are linearly independent over
$\Gcon$, otherwise we could replace one of them by a linear combination
of the others and thus decrease $n$. 

Embed $L$ in $K^{\imm}$, put $\Gimm = W(K^{\imm}) \otimes_{W(k)} \calO$,
and use Witt vector functoriality to construct
an embedding $\Galg \to \Gimm$. We can define the partial valuations
$v_n$, the valuations $w_r$
and the subring $\Gimmcon$ using the same formulas as before;
we then have $\Galgcon = \Galg \cap \Gimmcon$.

Let $w \in \Gperfcon$
be the Teichm\"uller lift of $t$. We show that 
every element $x$
of $\Gimm$ can be represented uniquely by a sum $\sum_{\alpha \in [0,1)}
x_{\alpha} w^{\alpha}$ with $x_{\alpha} \in \Gamma$ for all $\alpha$.
Namely, each element of $K^{\imm}$ can be written uniquely
as $\sum_{\alpha \in [0,1)} c_\alpha t^\alpha$ with
$c_\alpha \in k((t))$ for each $\alpha$. Thus we can choose
$x_{\alpha} \in \Gamma$ so that $\sum_{\alpha \in [0,1)} x_{\alpha} w^{\alpha}
\equiv x \pmod{\pi}$. But by the same reasoning applied to
$(x - \sum x_\alpha w_\alpha)/\pi$, we can choose the $x_{\alpha}$
so that $\sum x_\alpha w^{\alpha} \equiv x \pmod{\pi^2}$, and analogously
we can achieve the same congruence modulo any power of $\pi$. The limiting
values of $x_{\alpha}$ give the desired decomposition. This verifies
existence; as for uniqueness, it suffices to note that if $x \equiv 0
\pmod{\pi}$, by the uniqueness of the decompositions modulo $\pi$
we have $x_\alpha \equiv 0 \pmod{\pi^m}$ for all $m$, by induction on
$m$.

As above, decompose $x \in \Gimmcon$ as $\sum_{\alpha \in [0,1)}
x_{\alpha} w^{\alpha}$ with $x_{\alpha} \in \Gamma$; we claim that
in fact $x_\alpha \in \Gcon$ for all $\alpha$.
More specifically, choose $r,s>0$ such that $w_r(u/w) = 0$
and $w_r(x) \geq -s$; then we claim $r v_n(x_\alpha w^\alpha) + s + n
\geq 0$ for all $n$. Suppose this is not the case; choose the smallest $n$
for which $r v_n(x_\beta w^\beta) + s + n < 0$ for some $\beta$. For any
such $\beta$, we must have $v_p(x_{\beta}) = n$.
Write $x_\alpha = \sum_{i \in \ZZ} x_{\alpha,i} u^i$.
Let $y_\alpha$ be the sum of $x_{\alpha,i}u^i$ over all indices $i$ for
which $v_p(x_{\alpha,i}) < n$, and put $z_\alpha = x_\alpha - y_\alpha$
and $z = \sum_{\alpha} z_\alpha w^\alpha$.
Then $w_s(y_{\alpha}) \geq 0$ for each $\alpha$, so we also have
$w_r(z) \geq -s$. Also, $v_p(z) \geq n$ and $r v_n(z_\beta w^\beta)
+ s + n < 0$ for some $\beta$. But if we put $n = v_p(\pi^N)$, we then
have
\begin{align*}
rv_n(z_{\alpha+i} w^{\alpha+i}) + n &\geq
r v_0((z_{\alpha+i}/\pi^N) w^{\alpha+i}) + n \\
&= rv_0(z/\pi^N) + n \\
&\geq w_r(z) \geq -s
\end{align*} 
for any $\alpha \in [0,1)$ and $i \in \ZZ$. By choosing these so
$\alpha+i = \beta$, we obtain a contradiction.
We conclude that $w_{\alpha} \in \Gimmcon$ for each $\alpha \in [0,1)$.

Now write each $f_i$ as
$\sum_{\alpha \in [0,1)} f_{i,\alpha} u^\alpha$ with $f_{i,\alpha} \in \Gcon$;
then we have $\sum_{i=1}^n 
\sum_{\alpha \in [0,1)} f_{i,\alpha} g_i u^\alpha = 0$. By the uniqueness
of this type of representation, we conclude $\sum_{i=1}^n f_{i,\alpha}
g_i = 0$ for each $\alpha$. Since the $g_i$ were assumed to be linearly
independent over $\Gcon$, this implies $f_{i,\alpha} =0$ for each $i$
and $\alpha$, and so all of the $f_i$ are zero. This contradiction completes
the proof.
\end{proof}

The following lemma corresponds to de~Jong \cite[Corollary~8.2]{bib:dej1},
with essentially the same proof.
\begin{lemma}[after de~Jong] \label{lem:slope}
Let $M$ be a nonzero $\sigma$-module over $\Gcon = \GKcon$ with $K=k((t))$ and
$k$ algebraically closed, and let
$\phi: M \to \Gamma$ be a $\Gcon$-linear injective map such that
for some nonnegative integer $\ell$,
$\phi(F\bv) = p^{\ell} \phi(\bv)^\sigma$ for all $\bv \in M$. Then
the largest generic slope of $M$ is equal to $\ell$, occurring with
multiplicity~$1$. Moreover, $\phi^{-1}(\Gcon)$ is a $\sigma$-submodule of $M$
of dimension~$1$ with slope $\ell$.
\end{lemma}
\begin{proof}
By the previous lemma, the map $\phi: M \otimes_{\Gcon} \Galgcon
\to \Gamma \otimes_{\Gcon} \Galgcon \to \Galg$ is the composition
of two injections, so is injective. Let $s$ be the largest (generic)
slope of $M$ and $m$ its multiplicity. Choose $\lambda$ in a finite
extension of $\calO$ so that $\lambda$ is fixed by $\sigma$ and
$v_p(\lambda) = s$. By Proposition~\ref{prop:descfilt},
$M \otimes_{\Gcon} \Galgcon$ contains $m$
linearly independent eigenvectors $\bv_1, \dots, \bv_m$ of eigenvalue 
$\lambda$. Now
\[
\lambda \phi(\bv_i) = \phi(\lambda \bv_i) =
\phi(F\bv_i) = p^{\ell} \phi(\bv_i)^\sigma.
\]
This is only possible if $v_p(\lambda) = v_p(p^{\ell})$, i.e., if $s = \ell$.
In that case, we can take $\lambda = p^\ell$, in which case we must have
$\phi(\bv_i) \in \calO_0$ for each $i$; in particular, no two of the $\bv_i$
can be linearly independent. This is impossible unless $m=1$.

To complete the proof, it suffices to show that $\phi^{-1}(\Gcon)$ is
nonempty. If $F\bv = p^{\ell} \bv$ for $\bv \in M \otimes_{\Gcon} \Galgcon$,
we may choose a basis $\be_1, \dots, \be_n$ of $M$, write
$\bv = \sum c_i \be_i$, let $w \in \Gperfcon$ be the Teichm\"uller lift
of $t$, then imitate the proof of the previous lemma to write
each $c_i$ as a generalized power series
$\sum_{\alpha \in [0,1)} c_{i,\alpha} w^\alpha$ with $c_{i,\alpha} \in \Gcon$.
On one hand, $\phi(\bv) \in \calO$. On the other hand,
for each $\alpha$,
\[
\phi\left(\sum_i \sum_{\alpha \in [0,1)} c_{i,\alpha} w^{\alpha} \be_i \right) =
\sum_{\alpha \in [0,1)} w^{\alpha} \phi\left( \sum_i c_{i,\alpha} \be_i \right);
\]
since representations in the form $\sum_{\alpha \in [0,1)} d_\alpha w^\alpha$
are unique, we have $\sum_i c_{i,0} \be_i \in \phi^{-1}(\calO)$. Thus
$\phi^{-1}(\Gcon)$ is nonempty, and the proof is complete.
\end{proof}

\section{Local full faithfulness}

We now prove the local full faithfulness theorem, following
\cite[Theorem~9.1]{bib:dej1}. This theorem affirms 
\cite[Conjecture~2.3.2]{bib:tsu4}.
\begin{theorem} \label{thm:local}
Let $M_1$ and $M_2$ be $(\sigma, \nabla)$-modules over $\Gcon$,
and let $f: M_1 \otimes_{\Gcon} \Gamma \to 
M_2 \otimes_{\Gcon} \Gamma$ be a morphism of $(\sigma, \nabla)$-modules
over $\Gamma$. Then there exists a morphism $g: M_1 \to M_2$ that induces
$f$.
\end{theorem}
\begin{proof}
Regard $\nabla$ as a map from $M$ to itself by identifying
$\bv \in M$ with $\bv \otimes \frac{du}{u}$.
For $\ell$ sufficiently large, the Tate twist $M_1^*(\ell)$ of the dual
of $M_1$ is a $(\sigma, \nabla)$-module over $\Gcon$, and there is
a canonical isomorphism
$\Hom(M_1, M_2)(\ell) \cong M_1^*(\ell) \otimes M_2$.
Put $M = M_1^*(\ell) \otimes M_2$; then
$f$ induces an additive, $\Gcon$-linear map
$\phi: M \to \Gamma$ such
that:
\begin{enumerate}
\item[(a)] for all $\bv \in M$, $\phi(F\bv) = p^{\ell} \phi(\bv)$;
\item[(b)] for all $\bv \in M$, $\phi(\nabla \bv) = u \frac{d}{du} \phi(\bv)$.
\end{enumerate}
To prove the desired result, it suffices to prove that $\phi$ is induced
from a map $M \to \Gcon$ satisfying the analogues of (a) and (b),
i.e., that $\phi(M) \subseteq \Gcon$.
At this point, we may assume without loss of generality that
$k$ is algebraically closed.

Let $N \subseteq M$ be the kernel of $\phi$ on $M$; then $N$ is clearly
closed under $\sigma$ and $\nabla$ and saturated, so we may form the quotient
$(\sigma, \nabla)$-module $M/N$, and the induced map $\psi: M/N \to \Gamma$
is injective.
By Lemma~\ref{lem:slope},
the largest slope of $M/N$ is $\ell$ occurring
with multiplicity 1, and $P = \psi^{-1}(\Gcon)$ is
a sub-$\sigma$-module of $M/N$ of dimension 1 with slope $\ell$.
We now show that $P$ is also closed under $\nabla$.
If $\bv = \psi^{-1}(1)$,
then $F\bv = p^{\ell} \bv$ implies $\frac{du^\sigma/du}{u^\sigma/u}
F(\nabla \bv) = 
p^{\ell} \nabla\bv$, and so
\[
\frac{du^\sigma/du}{u^\sigma/u} p^{\ell} \psi(\nabla \bv)^\sigma
= p^{\ell} \psi(\nabla \bv).
\]
However, $\frac{du^\sigma}{du}$ is divisible by $\pi$ because
$u^\sigma \equiv u^p \pmod{\pi}$, so the two sides of the above
equation have different $p$-adic valuation unless $\psi(\nabla \bv) = 0$.
Since $\psi$ is injective, we conclude $\nabla \bv = 0$, so $P$ is
a $(\sigma, \nabla)$-submodule of $M/N$.

Since the slope $\ell$ of $P$ is greater than all of the other slopes of $M/N$,
we have an exact sequence
\[
0 \to P \to M/N \to (M/N)/P \to 0
\]
of $(\sigma, \nabla)$-modules
satisfying the conditions of Proposition~\ref{prop:split}; we conclude
that $M/N$ splits as a direct sum $P \oplus Q$ of
$(\sigma, \nabla)$-modules. If $Q$ is nonzero, we may
apply Lemma~\ref{lem:slope} once again to it, to conclude that its largest
slope is $\ell$, but this contradicts the fact that all slopes of $Q$
are smaller than $\ell$. Thus $Q$ must be the zero module and
$\phi(M) = \psi(M/N) = \psi(P) = \Gcon$, proving the desired result.
\end{proof}

\section{Rigid analytic Quillen-Suslin}

In the next section, we will reduce the global full faithfulness
theorem to a computation involving a finite projective module over the
ring $K\langle x_1,\dots, x_n \rangle^\dagger$ of
overconvergent power series in $n$ variables over $K$.
One might expect, in analogy to the Quillen-Suslin theorem, that such a module
must necessarily be free; since it is not too difficult to show that this
is actually the case, we include a proof here following (and ultimately
reducing to) the proof of the
Quillen-Suslin theorem given by Lang \cite{bib:lang}.

For an algebra $A$ complete with respect to a nonarchimedean absolute value
$|\cdot|$ and $\rho = (\rho_1, \dots, \rho_n)$ an $n$-tuple of positive reals,
we define $A \langle t_1, \dots, t_n \rangle_\rho$ as the ring
of formal power series which converge on the closed polydisc of radius
$\rho$; that is,
\[
A \langle t_1, \dots, t_n \rangle_\rho =
\left\{ \sum_I c_I t^I: c_I \in A, \lim_{\sum I \to \infty} |c_I| 
\rho^{I} = 0\right\}.
\]
Here $I = (i_1, \dots, i_n)$ runs over tuples of nonnegative integers,
$t^I = t_1^{i_1}\cdots t_n^{i_n}$,
$\rho^I = \rho_1^{i_1}\cdots \rho_n^{i_n}$,
and $\sum I = i_1 + \cdots + i_n$.
The ring $A \langle t_1, \dots, t_n \rangle_\rho$
is complete for the nonarchimedean absolute value
\[
\left| \sum_I c_I t^I \right| = \max_I \{ |c_I| \rho^{I} \}.
\]
For $n=1$, we define the \emph{leading term} $L(f)$ of a nonzero element
$f = \sum_{i=1}^\infty c_i t^i$ of $A \langle t \rangle_\rho$ as the
monomial $c_j t^j$ for $j$ the largest integer such that
$|c_j| \rho^{j} = \max_i \{|c_i|\rho^{i}\}$; we refer to $c_j$ as the
\emph{leading coefficient} of $f$. We define the \emph{degree} of $f$,
denoted $\deg(f)$, as the degree of its leading term.
Note that $|L(fg) - L(f) L(g)|
< |L(fg)|$; in particular, $\deg(fg) = \deg(f) + \deg(g)$.

We use two key lemmas to reduce from power series to polynomials. The
first is a form of Weierstrass preparation.
\begin{lemma}[Weierstrass preparation] \label{lem:weier}
Let $\gotho$ be a ring which is complete with respect to a nonarchimedean
absolute value $|\cdot|$, and put $A = \gotho \langle t \rangle_\rho$ for
some $\rho$.
Suppose the leading coefficient of $f \in A$ is a unit in $\gotho$. Then
there exists a unit $u$ in $A$ such that 
$u^{-1}f = \sum_{i=0}^j b_i t^i$ with
$b_j$ a unit and and $|b_i| \rho^{i} \leq |b_j| \rho^{j}$ for $i<j$
\end{lemma}
\begin{proof}
Put $f = \sum_{i \geq 0} f_i t^i$, $n = \deg(f)$ and
$g = \sum_{i \leq n} f_i t^i$.
Let $B$ be the ring of Laurent series $\sum_{i \in \ZZ} c_i t^i$ with
$c_i \in \gotho$ such that $|c_i| \rho^{i}$ remains bounded as $i \to
-\infty$ and goes to infinity as $i \to + \infty$.
Then as in 
\cite[Lemma~6.3]{bib:me7} (or \cite[Lemma~4.1.1]{bib:methesis}), $g^{-1} f$
factors in $B$ as $uv$ with $u = \sum_{i\geq 0} u_i t^i \in A$,
$v = 1 + \sum_{i < 0} v_i t^i$, $u_0$ a unit, and
$|v_i| \rho^{i} \leq 1$ for all $i$.
(Namely, one constructs a sequence of approximate factorizations which
converge under $|\cdot|$. The necessary initial condition is that
the leading term of $g^{-1} f$ be \emph{strictly} larger than all
other terms, which it is.)
In the equation $u^{-1} f = g v$, the left side
belongs to $A$ while the right side has no powers of $t$ beyond $t^n$.
Thus $u$ has the desired property.
\end{proof}
The second ingredient in our reduction is an argument commonly seen
in the proof of Noether normalization over finite fields (see
for instance \cite[Theorem~VIII.2.1]{bib:lang}). We state it two ways,
once for $\gotho$ a field, once for $\gotho$ not a field.
\begin{lemma} \label{lem:noeth1}
Let $\gotho$ be a complete discrete valuation field.
Let $\rho = (\rho_1, \dots, \rho_n)$ be a tuple of positive reals and
put $\rho' = (\rho_1, \dots, \rho_{n-1})$.
Suppose $u \in \gotho$ is a unit with $|u| \rho_n^{m} = 1$.
Put $B = \gotho \langle t_1, \dots, t_{n-1}
\rangle_{\rho'}$ and $A = \gotho \langle t_1, \dots, t_n \rangle_\rho
= B \langle t_n \rangle_{\rho_n}$.
Let $T_j: A \to A$ 
be the continuous $\gotho$-algebra isomorphism 
with 
\[
T_j(t_n) = t_n, \qquad T_j(t_i) = t_i + (ut_n^m)^{j^{n-i}}
\quad (i=1, \dots, n-1).
\]
Given $a \in A$, for all sufficiently large $j$ (depending on $a$),
the leading coefficient (in $t_n$) of $T_j(a)$ as an element of
$(\gotho \langle t_1, \dots, t_{n-1}
\rangle_{\rho'}) \langle t_n \rangle_{\rho_n}$
is a unit in $\gotho$.
\end{lemma}
\begin{proof}
Write $a = \sum_I a_I t^I$.
For $j$ sufficiently large, the leading
term in $T_j(a)$ will have degree $m (i_1 j^{n-1} + \cdots + i_{n-1} j^{n-1}
+ i_n)$, where $I = (i_1, \dots, i_n)$ is the last tuple in lexicographic
order that minimizes $|a_I| \rho^{I}$. Moreover, the coefficient will
be a power of $u$ plus smaller terms, which is a unit.
\end{proof}
\begin{lemma} \label{lem:noeth2}
Let $\gotho$ be a complete discrete valuation ring.
Let $\rho = (\rho_1, \dots, \rho_n)$ be a tuple of real numbers greater
than $1$, and
put $\rho' = (\rho_1, \dots, \rho_{n-1})$.
Put $B = \gotho \langle t_1, \dots, t_{n-1}
\rangle_{\rho'}$ and $A = \gotho \langle t_1, \dots, t_n \rangle_\rho
= B \langle t_n \rangle_{\rho_n}$.
Let $T_j: A \to A$ 
be the continuous $\gotho$-algebra isomorphism 
with 
\[
T_j(t_n) = t_n, \qquad T_j(t_i) = t_i + t_n^{mj^{n-i}}
\quad (i=1, \dots, n-1).
\]
Given $a \in A$, for all sufficiently large $j$ (depending on $a$)
and all sufficiently small $\lambda > 0$ (depending on $a$ and $j$),
the leading coefficient (in $t_n$) of $T_j(a)$ as an element of
$(\gotho \langle t_1, \dots, t_{n-1}
\rangle_{\lambda \rho'}) \langle t_n \rangle_{\rho_n^\lambda}$,
is a unit in $\gotho$.
\end{lemma}
\begin{proof}
This time, modulo $\idm$, the leading term of $T_j(a)$ (as a polynomial
in $t_n$) has unit coefficient; the same will be true of the leading term
of $T_j(a)$ within 
$(\gotho \langle t_1, \dots, t_{n-1}
\rangle_{\lambda \rho'}) \langle t_n \rangle_{\rho_n^\lambda}$
provided that $\lambda$ is sufficiently small.
\end{proof}
As a consequence of the dichotomy between Lemmas~\ref{lem:noeth1}
and~\ref{lem:noeth2}, the results we are about
to prove have two forms: one over a complete discrete valuation field,
in which $\rho$ does not change during the proof; and another over
a complete discrete valuation ring, in which the conclusion holds after
replacing $\rho$ by $\rho^\lambda$ for some $\lambda$ with $0 < \lambda \leq 
1$. For simplicity, we will state only the field versions explicitly
until we reach Theorem~\ref{thm:qs}.

An $n$-tuple $(f_1, \dots, f_n)$ of elements of a ring $R$
is \emph{unimodular} if its elements generate the unit ideal of $R$. 
Identifying $n$-tuples with column vectors, we say that two tuples
$f$ and $g$ are equivalent, notated $f \sim g$, if there exists
an invertible $n \times n$ matrix $M$ over $R$ such that $Mf = g$; this is
clearly an equivalence relation.
In this terminology, our analogue
of the main theorem of Quillen and Suslin is the following.

\begin{prop} \label{prop:unimod}
Let $K$ be a field complete with respect to a (nontrivial)
nonarchimedean absolute value
and let $\rho = (\rho_1, \dots, \rho_n)$ be a tuple such that for $i=1,
\dots, n$, some power of $\rho_i$ is the norm of an element of $K$.
Let $f$ be a unimodular tuple over $A = K \langle t_1, \dots, t_n 
\rangle_\rho$. Then $f \sim (1,0,\dots,0)$.
\end{prop}
We set some common notation for this proof and the next:
for $\rho = (\rho_1, \dots, \rho_n)$, put $\rho' = (\rho_1, \dots,
\rho_{n-1})$ and $B = K \langle t_1, \dots, t_{n-1} \rangle_{\rho'}$.
Also define $T_j$ as in Lemma~\ref{lem:noeth1} (or Lemma~\ref{lem:noeth2}
in the non-field case).
\begin{proof}
We prove the theorem by induction on $n$. If $n=0$, there is nothing to
prove. Suppose $n>0$; by Lemmas~\ref{lem:weier}
and~\ref{lem:noeth1}, for $j$ sufficiently large, 
each element of $T_j(f)$ is a unit in $A$ times
a polynomial over $B$ whose leading coefficient is a unit. By 
\cite[Theorem~XXI.3.4]{bib:lang}, $T_j(f)$ is equivalent over $A$
to a unimodular tuple over $B$. By the induction hypothesis, the latter
is equivalent to $(1,0,\dots,0)$. Hence $T_j(f) \sim (1,0,\dots,0)$; applying
$T_j^{-1}$ to the resulting matrix yields $f \sim (1,0,\dots,0)$, as desired.
\end{proof}

To apply this result to projective modules, we need an analogue of an
older result of Serre, which amounts to
the computation of $K_0$ of a polynomial ring.
(The hypothesis on $\rho$ ensures that
$K \langle t_1, \dots, t_n \rangle_\rho$ is noetherian;
see for instance \cite{bib:vdp}.)
\begin{prop} \label{prop:freeres}
Let $K$ be a field complete with respect to a (nontrivial)
nonarchimedean absolute value
and let $\rho = (\rho_1, \dots, \rho_n)$ be a tuple such that for $i=1,
\dots, n$, some power of $\rho_i$ is the norm of an element of $K$.
Then every finite module over $A = K \langle t_1, \dots, t_n \rangle_\rho$
has a finite free resolution.
\end{prop}
\begin{proof}
We proceed by induction on $n$; again, there is nothing to prove if $n=0$.
Note that by \cite[Theorem~XXI.2.7]{bib:lang},
if
\[
0 \to M_1 \to M \to M_2 \to 0
\]
is a short exact sequence of modules over a ring $R$ and any two
of $M, M_1, M_2$ have finite free resolutions, then so does the third.
Thus it suffices to show that every ideal $I$ of $A$ 
has a finite free resolution.

Let $f_1, \dots, f_m$ be generators of $I$.
By Lemmas~\ref{lem:weier} and~\ref{lem:noeth1}, for $j$ sufficiently large,
$T_j(f_i)$ can be written as a unit of
$A$ times a polynomial in $t_n$ over $B$.
Thus as an $A$-module, $I$ is isomorphic to $J \otimes_{B[t]} A$
for some ideal $J$ of $B[t]$. By \cite[Theorem~XXI.2.8]{bib:lang}
and the induction hypothesis that every finite module over $B$
has a finite free resolution, we deduce that $J$ has a finite free
resolution over $B[t]$. Hence $I$ has a finite free resolution over $A$,
completing the induction.
\end{proof}

Putting everything together, we get a result that implies the main
theorem of this section.
Note that the complex analogue of Proposition~\ref{prop:qs}
is also true; this follows
from a theorem of Lin \cite{bib:lin}.
\begin{prop} \label{prop:qs}
Let $K$ be a field complete with respect to a (nontrivial)
nonarchimedean absolute value and let $\rho = (\rho_1, \dots, \rho_n)$
be a tuple such that some power of each $\rho_i$ is the norm of an element
of $K$. Then every finite projective module over
$K \langle t_1, \dots, t_n \rangle_\rho$ is free.
\end{prop}
\begin{proof}
By Proposition~\ref{prop:freeres}, every finite projective module $M$
over
$K \langle t_1, \dots, t_n \rangle_\rho$ has a finite free resolution;
by \cite[Theorem~XXI.2.1]{bib:lang}, $M$ is stably free (the direct sum
of some finite free module with $M$ is finite free).
By \cite[Theorem~XXI.3.6]{bib:lang}, every stably free module over a 
ring with the unimodular extension property is free.
By Proposition~\ref{prop:unimod}, $K \langle t_1, \dots, t_n \rangle_\rho$
has the unimodular extension property. Thus every finite projective module
over $K \langle t_1, \dots, t_n \rangle_\rho$ is stably free and hence
free, as desired.
\end{proof}
\begin{theorem} \label{thm:qs}
Every finite projective module over $K \langle t_1, \dots, t_n \rangle^\dagger$
or $\calO \langle t_1, \dots, t_n \rangle^\dagger$ is free.
\end{theorem}
\begin{proof}
The ring $K \langle t_1, \dots, t_n \rangle^\dagger$ is the direct limit
of the rings $K \langle t_1, \dots, t_n \rangle_\rho$ over all tuples 
$\rho$ with $\rho_i > 1$ for all $i$, so
every finite projective over $K \langle t_1, \dots, t_n \rangle^\dagger$
is the base extension of a finite projective over $K \langle t_1, \dots,
t_n \rangle_\rho$ for some $\rho$. The result now follows from 
Proposition~\ref{prop:qs}.
For $\calO \langle t_1, \dots, t_n \rangle^\dagger$, the same argument holds
up to replacing $\rho$ by $\rho^\lambda$ for $0 < \lambda < 1$
unspecified (and replacing calls to Lemma~\ref{lem:noeth1} with
Lemma~\ref{lem:noeth2}), but the 
conclusion in the dagger algebra is unaffected.
\end{proof}

\section{Local to global}

We now proceed from the local full faithfulness theorem to a global statement,
by using a geometric lemma to ``push forward'' the problem from a general
variety to an affine space, where the reduction to the local theorem is 
straightforward.

For $X$ a smooth $k$-scheme of finite type and $\calE$ an overconvergent
$F$-isocrystal, let $H^0_F(X,\calE)$ denote the set of elements of
$\calE$ which are killed by $\nabla$ and fixed by $F$. Then
the full faithfulness of $j^*:
\FaIsocd(X/K) \to \FaIsoc(X/K)$ for $X$
follows from the fact that
the rank of $H^0_F(\calE, X)$
is the same whether computed in the overconvergent
or convergent category. The argument is the same as in Theorem~\ref{thm:local};
that is, given overconvergent $F$-isocrystals $\calE_1$ and $\calE_2$ on $X$,
$\Hom(\calE_1, \calE_2)$ is again an overconvergent $F$-isocrystal, and the
morphisms from $\calE_1$ to $\calE_2$ correspond to elements
of $\Hom(\calE_1, \calE_2)$ fixed by Frobenius and killed by $\nabla$,
in either the convergent or overconvergent category.

We now focus on showing that $\rank H^0_F(X,\calE)$ is the same in the convergent
and overconvergent categories.
The definition of $H^0_F(X,\calE)$ is local on $X$ in both cases, and by
the following lemma,
which follows from the main result of \cite{bib:etale2}, we may cover
$X$ with open affine subsets which are finite \'etale covers of affine spaces.
(The case of $k$ infinite and perfect is covered by \cite{bib:etale}; one
could reduce to this case with a bit of extra work.)
\begin{lemma}
Let $X$ be a separated $k$-scheme of finite type of pure dimension $n$
and let $x$ be a smooth
geometric point of $X$. Then there exists an open dense subset $U$ of $X$,
containing $x$ and defined over $k$, such that $U$ admits a finite \'etale
map over $k$ to affine $n$-space.
\end{lemma}
Thus it suffices to consider the case $X$ equal to such an open dense
subset $U$. Let $f: X \to \AAA^n$
be a finite \'etale morphism and $\calE$ an overconvergent $F$-isocrystal
on $X$.
There is a pushforward construction in the overconvergent
and convergent categories \cite[Proposition~5.1.2]{bib:tsu4}
such that $H^0(X,\calE) = H^0(\AAA^n, f_* \calE)$.
Thus it suffices to consider isocrystals on affine space itself.

Take $f = a$, so that $q = p^a$.
Let $R = \calO\langle x_1,\dots, x_n \rangle^\dagger$ be the ring of
overconvergent power series in $n$ variables over $\calO$,
and let $\sigma$
be the Frobenius lift on $R$ sending $\sum_I c_I x^I$ to $\sum_I 
c_I^\sigma x^{pI}$. Then the data of
an overconvergent $F$-isocrystal $\calE$
on $\AAA^n$ is simply that of
a $(\sigma, \nabla)$-module $M$ on $R[\fp]$; we can find an integer
$\ell$ such that the Tate twist $M(\ell)$ can be defined over $R$.
Then $M(\ell)$ is free over $R$ by Theorem~\ref{thm:qs}.
Let $R^\wedge$ denote the
$\pi$-adic completion of $R$; then the elements of $H^0_F(\AAA^n,\calE)$ in
the convergent and overconvergent categories correspond to
the elements $\bv \in M(\ell)
\otimes_{R} R^\wedge$ and $M(\ell)$, respectively,
such that $F\bv = p^{\ell} \bv$ and $\nabla \bv = 0$.

Note that a power series in $x_1, \dots, x_n$ over $\calO$ lies in $R$
if and only if it is overconvergent in each variable separately.
In other words, if $S_i$ is the valuation subring of
the fraction field of the ring of
null power series in all of the $x$'s other than $x_i$, and
$R_i = S_i\langle x \rangle^\dagger$, then $\cap_i R_i = R$
(the intersection taking place in the completed fraction field of $R$).

We use the above observation to ``take apart'' $M$. Namely,
by definition, the module $\Omega^1_{R}$ is freely generated over $R$
by $dx_1,\dots, dx_n$. Let $\Omega^1_{R_i/S_i}$ be the free
module generated over $R_i$ by $dx_i$. If we put $k = S_i$,
then $M(\ell) \otimes_{R} R_i$
has a natural structure as $(\sigma, \nabla_i)$-module over
$R_i$, with $\nabla_i$ given by composing the given map
$M(\ell) \to M(\ell) \otimes \Omega^1_R = \oplus_i R \,dx_i$ with the
projection onto the $i$-th factor and the map $R\,dx_i \to R_i\,dx_i$.

Identify $R_i$ with a subring of $\Gamma^{S_i}_{\con}$ so
that $x_i^{-1}$ reduces to a uniformizer of the residue field.
If $\bv \in M(\ell) \otimes_R R^\wedge$ satisfies $F\bv = p^\ell \bv$
and $\nabla \bv = 0$,
then Theorem~\ref{thm:local}
applied to $M(\ell) \otimes_R R_i$ implies that $\bv \in M(\ell)
\otimes_R R_i$ for each $i$.
Since $M(\ell)$ is free over $R$ and $\cap R_i = R$,
we deduce $\bv \in M(\ell)$.
We thus conclude that
$H^0(\AAA^n, \calE)$ has the same rank in the convergent and
overconvergent categories; as noted above, this suffices to prove
Theorem~\ref{thm:main}.

\end{document}